\documentclass[11pt]{amsart}
\usepackage{amsmath, amsfonts, amsthm, amssymb, multicol}
\usepackage{graphicx}
\usepackage{float}
\usepackage{verbatim}

\allowdisplaybreaks

\usepackage[square,sort,comma,numbers]{natbib}
\usepackage{fancyhdr}
 
\numberwithin{equation}{section}

\newtheorem{thm}{Theorem}[section]
\newtheorem{lma}[thm]{Lemma}
\newtheorem{cor}[thm]{Corollary}

\newtheorem{ques}[thm]{Question}

\renewcommand{\epsilon}{\varepsilon}
\newcommand{\eps}{\varepsilon}

\newcommand{\rd}{\mathbb{R}^d}

\renewcommand{\geq}{\geqslant}
\renewcommand{\leq}{\leqslant}

\newcommand{\ad}{\dim_{\mathrm{A}} }
\newcommand{\as}{\dim^\theta_{\mathrm{A}} }
\newcommand{\ls}{\dim^\theta_{\mathrm{L}} }

\newcommand{\ld}{\dim_{\mathrm{L}}  }

\newcommand{\mld}{\dim_{\mathrm{ML}}  }
\newcommand{\hd}{\dim_{\mathrm{H}}  }
\newcommand{\bd}{\dim_{\mathrm{B}}  }
\newcommand{\pd}{\dim_{\mathrm{P}}  }

\title{Fractal geometry of \\ Bedford-McMullen carpets}

\author{Jonathan M. Fraser\\ \\
 U\MakeLowercase{niversity of} S\MakeLowercase{t} A\MakeLowercase{ndrews}, S\MakeLowercase{cotland} \\
\MakeLowercase{Email: jmf32@st-andrews.ac.uk}}
\thanks{The  author was   supported by an \emph{EPSRC Standard Grant} (EP/R015104/1) and a  \emph{Leverhulme Trust Research Project Grant} (RPG-2019-034).}

\begin{document}


\maketitle
\thispagestyle{empty}

\begin{abstract}
In 1984 Bedford and McMullen independently introduced a family of self-affine sets now known as \emph{Bedford-McMullen carpets}.  Their work stimulated a lot of research in the areas of fractal geometry and non-conformal dynamics.  In this survey article we discuss some aspects of  Bedford-McMullen carpets, focusing mostly on dimension theory.  
\\ 


\end{abstract}

\section{Bedford-McMullen carpets}

One of the most important and well-studied methods for generating interesting fractal sets is via iterated function systems (IFSs).  Roughly speaking, an \emph{IFS} is a finite collection of contraction mappings acting on a common compact domain, and the associated \emph{attractor}  is the unique non-empty compact set which may be expressed as the union of   scaled down copies of itself under the maps in the IFS.  \emph{Self-similar sets} are attractors of IFSs where the contractions are similarities, and \emph{self-affine sets} are attractors of IFSs where the maps act on a Euclidean domain and are affine contractions  (the composition of a linear contraction and a translation).   See \cite{falconer} for more background on IFSs and the survey \cite{affinesurvey} for a detailed history of self-affine sets and measures as well as \cite{baranyinvent} for a recent breakthrough in the dimension theory of general self-affine sets.

Affine maps may scale by different amounts in different directions (as well as skewing and shearing) and this leads to self-affine sets being rather more complicated than self-similar sets.   Bedford-McMullen carpets are the simplest possible family of (genuinely) self-affine sets.  They preserve the key feature of self-affinity: different scaling in different directions, but everything else about the construction is as simple as possible.  The simplicity of the model, combined with  the ability to capture a key aspect of the theory, has contributed greatly to its popularity.  In fact, Bedford-McMullen carpets provide an excellent example to aspiring mathematicians: a good model should reveal and capture a new phenomenon, but remain as simple as possible.

\begin{figure}[h]
	\centering
	\includegraphics[width=\textwidth]{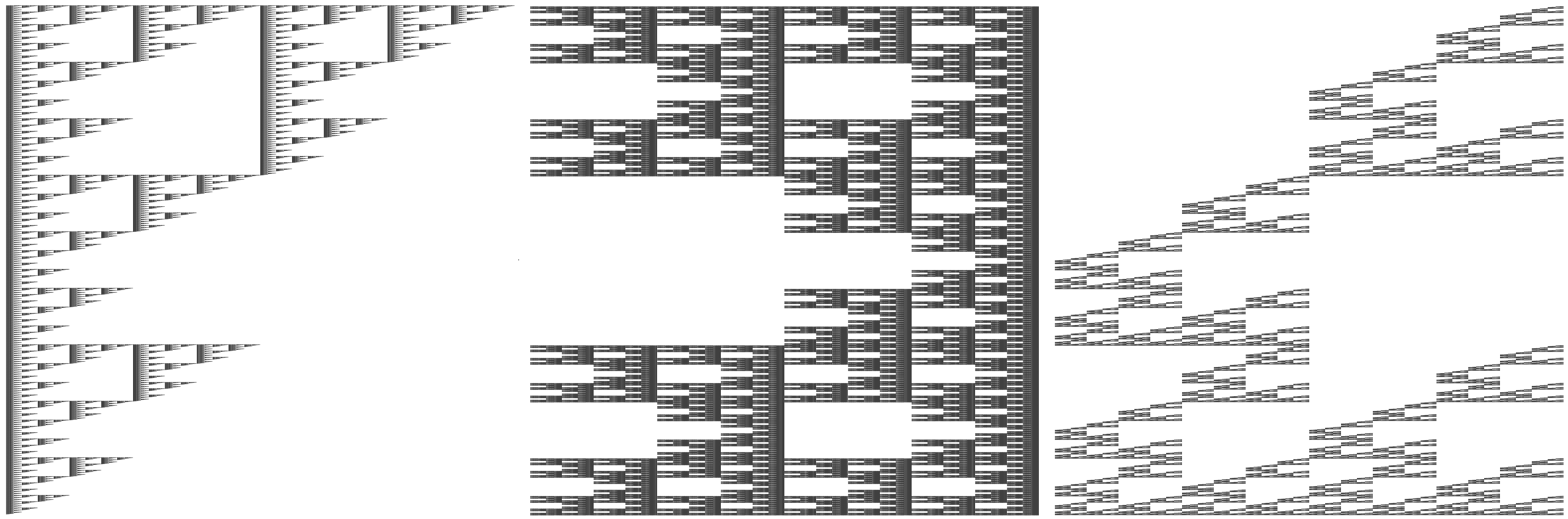}
\caption{Three examples of  Bedford-McMullen carpets based on the $2 \times 3$ grid.  } \label{bmc1}
\end{figure}

Let us first  recall the  Bedford-McMullen construction following \cite{bedford, mcmullen}. We work in the Euclidean plane, and begin with the unit square, $[0,1]^2$.  Fix integers  $n>m>1$, and divide  the unit square   into an $m\times n$ grid.  Select a subset of the rectangles formed by the grid and consider the IFS consisting of the  affine maps which map $[0,1]^2$ onto each chosen rectangle preserving orientation (that is, the affine part of each map is the diagonal matrix with diagonal entries $1/m$ and $1/n$).  The attractor of this IFS is a self-affine set, and such self-affine sets are known as  \emph{Bedford-McMullen carpets}, see Figure \ref{bmc1}. 

Bedford-McMullen carpets also have an important role in the theory of expanding dynamical systems.  Viewed as subsets of the 2-torus $[0,1)^2$, Bedford-McMullen carpets are invariant under the toral endomorphism $(x,y) \mapsto (m x \text{ mod } 1, \, n y \text{ mod } 1)$, which provides  a simple model of a non-conformal dynamical system.   Since the work of Bedford and McMullen, the study of self-affine carpets has  received sustained interest in the literature. Generally speaking a `carpet' is an attractor of an IFS acting on the plane consisting of affine maps whose linear parts are given by diagonal matrices (or possibly anti-diagonal matrices).  There are now many popular  families of carpet, generalising the Bedford-McMullen model in various ways.  Lalley-Gatzouras carpets \cite{lalley-gatz} maintain the column structure but allow the diagonal matrix to vary, Bara\'nski carpets  \cite{baranski} maintain the grid structure but allow the matrices to vary.  A crucial difference between the Bara\'nski  and Lalley-Gatzouras models is that Bara\'nski allows the strongest contraction to be in either direction, whereas Lalley-Gatzouras insists that the strongest contraction be in the vertical direction.  Generalising both Lalley-Gatzouras and Bara\'nski carpets is the family introduced by Feng-Wang \cite{fengwang} which allows arbitrary non-negative diagonal matrices and generalising the Feng-Wang family is a family we introduced which allows arbitrary diagonal and anti-diagonal matrices \cite{boxlike}.   There are also models which step out of the carpet programme whilst maintaining several of the key features, such as excessive alignment of cylinders, for example, \cite{istvan2}. In order to keep this survey concise, we  make no further mention of carpets outside the Bedford-McMullen family.  There is a vast literature on Bedford-McMullen carpets and as such there are a lot of interesting research directions which we will not discuss here.  This survey is mostly focused on the fractal geometry and dimension theory of Bedford-McMullen carpets and associated self-affine measures.

\section{Dimension theory}

A central aspect of fractal geometry is dimension theory.  In fractal settings, fine structure makes the task of simply  \emph{defining} dimension already an interesting problem.  Roughly speaking, a `dimension' should describe how an object fills up space on small scales.  There are many ways to describe this, however, and an important aspect  of the subject is in understanding the relationships and differences between the many different notions of dimension; each of interest in its own right.    In this survey we   focus on Hausdorff, packing, box, Assouad and lower dimension, which we denote by $\hd, \pd, \bd, \ad, \ld$, respectively.  Often one needs to consider upper and lower box dimension separately, but for the sets we discuss these coincide and so we brush over this detail.  We will not define these notions here, but refer the reader to \cite{falconer, mattila} for the definitions and an in depth discussion of the Hausdorff, packing, and box dimensions and \cite{fraserbook} for the Assouad and lower dimensions.   It is useful to keep in mind that for all non-empty compact  sets $E \subseteq \rd$ (with equal upper and lower box dimension),
\[
0 \leq \ld E \leq \hd E \leq \pd E \leq \bd E \leq \ad E \leq d.
\]
Bedford \cite{bedford} and McMullen \cite{mcmullen} independently obtained explicit formulae for the   Hausdorff, packing, and box dimensions of Bedford-McMullen carpets.  More recently, in 2011,  Mackay \cite{mackay} computed the Assouad dimension, and, in 2014, Fraser computed the lower dimension \cite{fraserassouad}.   We need more notation in order to state these results.  

Let $N$ be the number of maps in the defining IFS (that is the number of chosen rectangles), and let  $M$ be the number of columns containing at least one chosen rectangle.  Finally, let   $N_i>0 $ be the number of rectangles  chosen from the $i$th  non-empty column. 

\begin{thm} \label{bedmcformulae}
Let $F$ be a Bedford-McMullen carpet.  Then
\begin{align*}
\ad F  \ &  = \  \frac{\log M}{\log m} \, + \, \max_{i} \frac{\log N_i}{\log n}, \\ \\
 \pd F \ = \ \bd F   \ &= \  \frac{\log M}{\log m} \, + \, \frac{\log ( N/M)}{\log n}, \\ \\
\hd F  \ &= \ \frac{\log \sum_{i=1}^{M} N_i^{\log m/ \log n}}{\log m}  \\ \\
\text{and} \qquad \qquad \ld F  \ &= \  \frac{\log M}{\log m} \, + \, \min_{i} \frac{\log N_i}{\log n}.
\end{align*}
\end{thm}

\emph{Sketch proof.}  We sketch the argument giving the box dimension and only discuss the rough ideas for the others.  The box dimension of a bounded set $E$ captures the polynomial  growth rate of $N_r(E)$ as $r \to 0$, where $N_r(E)$ denotes the smallest number of open sets of small diameter $r \in (0,1)$ required to cover $E$.  That is, the box dimension can loosely be defined by $N_r(E) \approx r^{-\bd E}$.

Let $r>0$ be very small and  $k$ be an integer such that $r \approx n^{-k}$.  The $k$th level cylinders in the construction of the carpet  $F$ are rectangular sets of height $\approx r$ and length $m^{-k}$ (which is rather longer). Therefore, when  looking for optimal $r$-covers of $F$ we may treat the $k$th level cylinders separately. Let $l$ be an integer such that $r \approx m^{-l}$ and consider the $l$th level cylinders inside a given $k$th level cylinder.  This forms a grid and cylinders in the same column may be covered efficiently by $\approx 1$ set of diameter $r$ and therefore we only need to count the non-empty columns, of which there are $M^{l-k}$.  The total number of $k$th level cylinders is $N^k$ and therefore
\[
N_r(F) \approx N^k M^{l-k} \approx r^{-(\log N/\log n+\log M/\log m -\log M/\log n)}
\]
as required. 

The Hausdorff dimension is more awkward to compute. The lower bound is usually handled via measures, either by the mass distribution principle (see \cite[Chapter 4]{falconer}) or by direct computation of the Hausdorff dimension of a suitable measure.  In fact we sketch this part of the proof in Section \ref{secmeas}.  The upper bound is proved by a delicate covering argument.  The key difference between this argument and the covering argument we sketched for box dimension above is that Hausdorff dimension allows different sized sets in the cover and so the difficulty is in deciding which cylinders to cover together and which to break up into smaller pieces.  A more direct approach to proving the upper bound, which ultimately boils down to constructing a delicate cover, is to show that the lower local dimension of a suitable measure  is at most $h$ for all points $x \in F$, where $h=\hd F$ is the intended Hausdorff dimension.   In contrast to the mass distribution principle, which asks for the measure of a ball never to be too big, this approach asks for balls around all points to have large mass infinitely often.   See  \cite[Chapter 4]{falconer} for more on this approach to finding upper bounds for Hausdorff dimension in general. The McMullen measure, defined later by \eqref{mcmeas},  can be used for both the upper and lower bounds.

The Assouad and lower dimensions are in some sense dual to each other and so we only discuss Assouad dimension.  The lower bound is most efficiently proved via \emph{weak tangents}.  See \cite[Chapter 5]{fraserbook} for more on weak tangents in the context of dimension theory.     Mackay \cite{mackay} constructed a weak tangent which is the product of two self-similar sets, one of dimension $ \log M/\log m$ (the projection of $F$ onto the first coordinate) and the other of dimension  $\max_{i} \log N_i/\log n$ (the maximal vertical slice of $F$).  Then the lower bound follows since the Hausdorff dimension of a weak tangent is a lower bound for the Assouad dimension.  Mackay proved the upper bound by a direct covering argument, similar to that sketched above for box dimension.  An alternative argument giving the upper bound  using the Assouad dimension of measures is given in \cite{fraserhowroyd}. \hfill \qed 
\\

The tangent structure of self-affine sets is particularly interesting since the small scale structure typically differs greatly from the large scale structure.  This is very different from self-similar sets, for example.  As mentioned in the above proof, Mackay \cite{mackay} used tangent sets with a product structure to study the Assouad dimension.  This product structure is seen much more generally.  Bandt and K\"aenm\"aki \cite{bandtkaenmaki} gave  a general description of the tangent structure of Bedford-McMullen carpets in the case where $M=m$.  This result has been generalised in various ways, see  \cite{algomhochman, fibred, kaenmaki, rossi}. 

Returning to Theorem \ref{bedmcformulae}, the following amusing commonality was pointed out to me by Kenneth Falconer.  For $p \in [0, \infty] \cup \{-\infty\}$, write
\[
\|\underline{N} \|_p =   \left( \frac{1}{M} \sum_{i=1}^{M} N_i^p \right)^{1/p}
\]
for the ``$p$th average'' of the vector $\underline{N} : = (N_1, \dots, N_M)$ describing the number of rectangles in each of the non-empty columns. We adopt the natural interpretation of $\|\underline{N} \|_\infty = \max_i N_i$ and $\|\underline{N} \|_{-\infty} = \min_i N_i$. Then, the expression
\begin{equation} \label{pdim}
 \frac{\log M}{\log m} \, + \,  \frac{\log \|\underline{N} \|_p}{\log n}
\end{equation}
gives the Assouad dimension when $p= \infty$, the box and packing dimensions when $p=1$, the Hausdorff dimension when $p=\log m/\log n$, and the lower dimension when $p=-\infty$.  This observation warrants the question of whether there are sensible, perhaps yet undiscovered, notions of fractal dimension corresponding to other values of $p$.  Moreover, the expression  \eqref{pdim} has a useful interpretation as the dimension of the projection of $F$ onto the first coordinate plus the `average' column dimension.

It is no surprise that the box and packing dimensions coincide in Theorem \ref{bedmcformulae}.  Indeed, the packing dimension and upper box dimensions coincide much more generally: see \cite[Corollary 3.10]{falconer}, which applies to very general IFS attractors.  This identity aside, we see that Bedford-McMullen carpets provide an excellent model for understanding the differences between the different notions of dimension. If all non-empty columns contain the same number of rectangles  (that is, $N_i=N/M$ for all $i$), then we say the carpet has \emph{uniform fibres} and otherwise it has \emph{non-uniform fibres}.  There is a simple dichotomy:  in the uniform fibres case 
\[
\ld F = \hd F = \bd F = \ad F
\]
and, in the non-uniform fibres case, 
\[
\ld  F < \hd F < \bd F < \ad F .
\]
The question of Hausdorff and packing \emph{measure} for Bedford-McMullen carpets is subtle.  Peres \cite{perespacking, pereshausdorff} proved that in the non-uniform fibres case both the Hausdorff and packing measures are infinite in their respective dimensions.   It is instructive to compare this with the situation for self-similar sets where the open set condition is enough to ensure that the Hausdorff and packing measures are positive and finite in their dimension. We write  $\mathcal{P}^h$ and $\mathcal{H}^h$ for the packing and Hausdorff measures with respect to a gauge function $h$, see \cite{falconer}. In \cite{perespacking} it is shown that $\mathcal{P}^h(F)=\infty$ for
\[
h(x) = \frac{x^{\pd F}}{|\log x|}
\]
but $\mathcal{P}^h(F)=0$ for
\[
h(x) = \frac{x^{\pd F}}{|\log x|^{1+\varepsilon}}
\]
for all $\varepsilon \in (0,1)$.   In \cite{pereshausdorff} it is shown that $\mathcal{H}^h(F)=\infty$ for
\[
h(x) = x^{\hd F}\exp \left(\frac{-c |\log x|}{(\log |\log x|)^2}\right)
\]
with $c>0$ small enough, but $\mathcal{H}^h(F)=0$ for
\[
h(x) = x^{\hd F}\exp \left(\frac{- |\log x|}{(\log |\log x|)^{2-\varepsilon}}\right)
\]
for all $\varepsilon \in (0,1)$.  In particular these results show that
\[
\mathcal{H}^{\hd F} (F) = \mathcal{P}^{\pd F} (F) = \infty .
\]

\section{Interpolating between dimensions}

A new perspective in dimension theory is that of `dimension interpolation', see  \cite{frasersurvey}.  Roughly speaking, the idea is to consider two distinct notions of dimension $\dim$ and $\text{Dim}$, which satisfy $\dim E \leq \text{Dim} \, E$ for all sets $E \subseteq \rd$, and introduce a continuously parametrised family of dimensions $\dim_\theta$ for $\theta \in (0,1)$ which satisfy $\dim E \leq \dim_\theta E \leq  \text{Dim} \,  E$.  Crucially, the dimensions $\dim_\theta$ should capture some key features of both $\dim$ and $\text{Dim}$ in a geometrically interesting way and the hope is to provide more nuanced information than that provided by $\dim$ and $\text{Dim}$ when considered in isolation.    Since the most commonly studied dimensions are typically distinct in this case, Bedford-McMullen carpets provide the ideal testing ground for this approach.  

\subsection{The Assouad spectrum} 

The Assouad spectrum, introduced in \cite{fyu} and denoted by $\as$, interpolates between the (upper) box dimension and the (quasi-)Assouad dimension.  The Assouad dimension of a bounded set $E \subseteq \mathbb{R}^d$  is defined by considering $N_r(B(x,R) \cap E)$, that is, the size of an optimal $r$-cover  of an $R$-ball where $0<r<R$ are two independent scales. The Assouad spectrum fixes the relationship between these two scales by forcing $R=r^\theta$, where $\theta \in (0,1)$ is the interpolation parameter.  The result is a family of dimensions $\as E$ which is continuous in $\theta$ and satisfies $  \bd E \leq \as E \leq \ad E$ for all $\theta \in (0,1)$.  The analogous lower spectrum $\ls$ is defined by a similar modification of the definition of lower dimension.  The Assouad and lower spectra of Bedford-McMullen carpets were computed in \cite{fyu2}.  We write $N_{\max} = \max_i N_i$ and $N_{\min} = \min_i N_i$.

\begin{thm} \label{CarpetsSpectra}
Let $F$ be a Bedford-McMullen carpet.  Then, for all $0<\theta \leq \log m / \log n $,
\[
\dim_{\mathrm{A}}^\theta F \ = \    \frac{  \log M  -   \theta  \log ( N /N_{\max})  }{(1- \theta)\log m}  
+ 
\frac{   \log ( N/M) -   \theta  \log N_{\max}   }{(1- \theta)\log n} 
\]
and
\[
\dim_\mathrm{L}^\theta  F \ = \     \frac{  \log M  -   \theta  \log ( N /N_{\min})  }{(1- \theta)\log m}  
+ 
\frac{   \log ( N/M) -   \theta  \log N_{\min}   }{(1- \theta)\log n} 
\]
and,  for all $  \log m / \log n  \leq \theta < 1$,
\[
\dim_{\mathrm{A}}^\theta F \ = \    \frac{\log M}{\log m} \, + \,   \frac{\log N_{\max}}{\log n}
\]
and
\[
\dim_\mathrm{L}^\theta  F \ = \   \frac{\log M}{\log m} \, + \,   \frac{\log N_{\min}}{\log n}.
\]
\end{thm}

In both the Assouad and lower spectrum there is a single phase transition at $\theta = \frac{\log m}{\log n}$. In many  other examples where the Assouad spectrum is known there is a similar phase transition occurring at a value of $\theta$  with  a particular geometric significance, see \cite{fraserbook}.  In this case it is the ratio of the Lyapunov exponents (for any ergodic measure) or the `logarithmic eccentricity' of the cylinders in the construction.

\begin{figure}[h]
	\centering
	\includegraphics[width=0.47\textwidth]{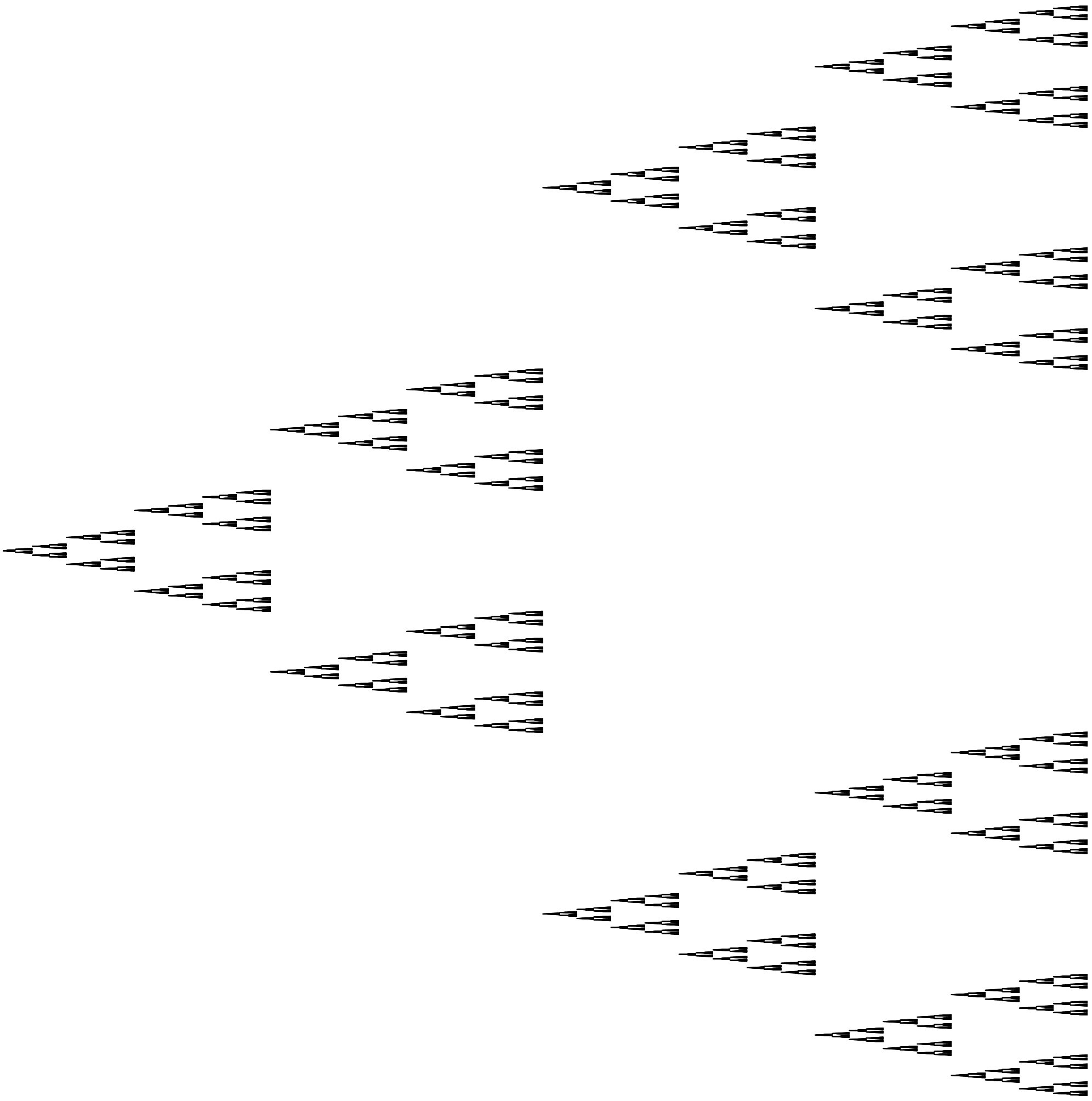} \quad
\includegraphics[width=0.47\textwidth]{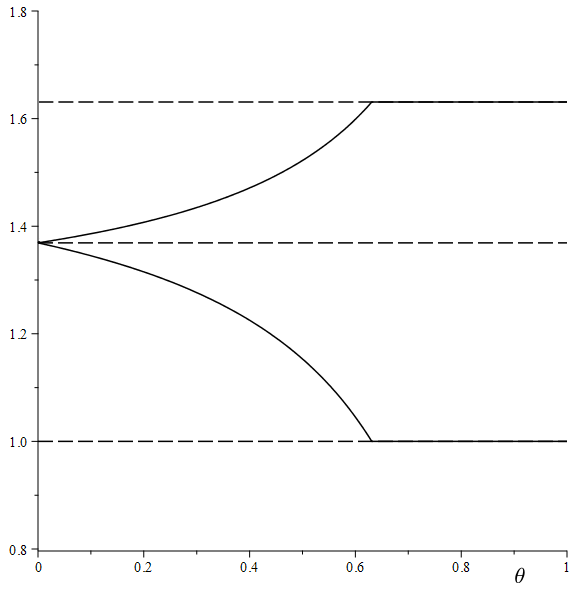} \\
\caption{Left: a  Bedford McMullen carpet with $m=2$,  $n=3$ and and $N_1=1<N_2=2$.  Right: plots of the Assouad and lower spectra.  Dotted lines at the lower, box, and Assouad dimensions are shown for comparison. }
\end{figure}

The curious reader may at this point wonder if the Assouad spectrum is a partial solution to the problem of finding dimensions satisfying \eqref{pdim}.  By continuity and monotonicity, for each $0<\theta < \frac{\log m}{\log n}$ there is a unique $p \in (1, \infty)$ such that 
\[
\as F =  \frac{\log M}{\log m} \, + \,  \frac{\log \|\underline{N} \|_p}{\log n},
\]
however, the function mapping $\theta$ to $p$ is not satisfying: it depends too strongly on $F$.  Ideally we would want a function depending only on $\theta$ and possibly $m$ and $n$.

\subsection{Intermediate dimensions}

The intermediate dimensions, denoted by $\dim_\theta$ and introduced in \cite{FalconerFraserKempton}, interpolate between the (upper or lower) box dimension and the Hausdorff dimension.  The Hausdorff and box dimension are both defined by considering efficient covers of the set.  The Hausdorff dimension places no restriction on the relative sizes of the sets used in the cover and weights their contribution to the dimension according to their size (see definition of Hausdorff measure \cite{falconer}) and box dimension considers covers by sets of the same size.  The intermediate dimensions impose partial restrictions on the relative sizes of the covering sets by insisting that $|U| \leq |V|^\theta$ for all covering sets $U,V$ with diameters $|U|, |V| \leq 1$.    The intermediate dimensions are continuous in $\theta \in (0,1]$ with $\dim_1 E = \bd E$ and satisfy
\[
\hd E \leq \dim_\theta E \leq \bd E
\] 
for all bounded $E \subseteq \rd$.   The intermediate dimensions are not necessarily continuous at $\theta = 0$, that is, they do not necessarily approach the Hausdorff dimension as $\theta \to 0$.  Establishing continuity at 0 for particular examples turns out to be  a key problem.  For example, if the intermediate dimensions are continuous at  0, then strong applications can be derived concerning the box  dimensions of projections and images under stochastic processes, see \cite{burrell1, burrell2}.  

Computing an explicit formula for the intermediate dimensions of Bedford-McMullen carpets seems to be a difficult problem, investigated in \cite{FalconerFraserKempton} and \cite{istvan}. See also the survey \cite{falconersurvey}.  We summarise what is known so far, for $F$ a Bedford-McMullen carpet with non-uniform fibres:
\begin{enumerate}
\item $\dim_\theta F$ is continuous at 0, and therefore interpolates between the Hausdorff and box dimensions \cite{FalconerFraserKempton};
\item for all $\theta \in (0,1)$, $\dim_\theta F > \hd F$ \cite{FalconerFraserKempton};
\item for all $\theta \in (0,1)$, $\dim_\theta F < \bd F$ \cite{istvan};
\item there are upper and lower bounds for $\dim_\theta F$ with the best known given in \cite{istvan};
\item $\dim_\theta F$ is not necessarily concave \cite{istvan}, which is noteworthy since most of the basic examples have turned out to have concave intermediate dimensions.
\end{enumerate}

\section{Invariant measures} \label{secmeas}

The interplay between invariant sets and invariant measures is central in fractal geometry and ergodic theory.  There are many natural $(\times m, \times n)$ invariant measures living on Bedford-McMullen carpets but perhaps the most natural are the self-affine  measures.  Given a Bedford-McMullen carpet with first level rectangles indexed by the set $\{1, \dots, N\}$, associate a probability vector $\{p_1, \dots, p_N\}$, that is, with $0<p_d<1$ for all $d$ and $\sum_d p_d = 1$.  Let $\mu$ be the measure formed by iteratively subdividing unit measure among the $N$ rectangles at each stage in the construction of the carpet $F$ according to the probability vector.  More formally, let $\mathbb{P}$ denote the Bernoulli measure on the symbolic space $\{1, \dots, N\}^{\mathbb{N}}$ consisting of infinite one-sided words over the alphabet $\{1, \dots, N\}$.  Then $\mu = \mathbb{P} \circ \Pi^{-1}$ where $\Pi$ is the associated coding map which sends a word $(d_1, d_2, \dots)$ to the  point in $F$ coded by the corresponding sequence of rectangles.    In particular, $\mu$ is a Borel probability measure fully supported on $F$ and invariant under the endomorphism $(\times m, \times n)$.  The measure $\mu$ is a self-affine measure since it is the unique Borel probability measure invariant under the IFS weighted by the probability vector.  There is a rich dimension theory of measures, which we will not dwell on here.  Measures invariant under nice enough dynamical systems are often `exact dimensional', which means that many of the familiar notions of dimension for measures coincide (e.g. Hausdorff, packing, entropy dimensions).  A Borel measure $\nu$ is \emph{exact dimensional} (with dimension $\alpha$) if for $\nu$ almost all $x$
\begin{equation} \label{ldd}
\lim_{r \to 0} \frac{\log \nu (B(x,r))}{\log r} = \alpha.
\end{equation}
The expression on the left is the \emph{local dimension} of $\nu$ at $x$ (when the limit  exists) and so exact dimensionality can be characterised as the local dimension being almost surely constant.  Self-affine measures on Bedford-McMullen carpets are known to be exact dimensional and, moreover, the dimension satisfies the \emph{Ledrappier-Young formula}.  The Ledrappier-Young formula stems from influential papers \cite{led, led2} which established a deep connection between dimension, entropy and Lyapunov exponents in the general context of measures invariant under $C^1$-diffeomorphisms.  Recently there has been a lot of progress establishing Ledrappier-Young formulae for measures invariant under expanding maps, such as those associated with IFSs.  In particular, it is now known that \emph{all} self-affine measures are  exact dimensional and satisfy the appropriate Ledrappier-Young formula,  see \cite{barany, baranykaenmaki, feng}.   The case of self-affine measures on Bedford-McMullen carpets was resolved in \cite{kp-measures}.  In this setting the   Lyapunov exponents are $\log m< \log n$ and the entropy is given by
\[
h(\mu) =  -\sum_{d=1}^N p_d \log p_d.
\]
The entropy of the projection $\pi \mu$ of $\mu$ onto the first coordinate also plays a role.  Here we have to sum weights belonging to the same column and so the entropy is 
\[
h(\pi \mu) =  -\sum_{i=1}^M \Big(\sum_{d : \pi d = i} p_d\Big) \log\Big(\sum_{d : \pi d = i} p_d\Big).
\]

\begin{thm} \label{lyformula}
Self-affine measures $\mu$ on Bedford-McMullen carpets are exact dimensional with dimension given by
\[
\dim \mu =  \frac{h(\pi \mu)}{\log m} + \frac{h(\mu)-h(\pi \mu)}{\log n}.
\]
\end{thm}
We give a sketch proof here, which is deliberately not at all rigorous but hopefully shows where the formula comes from and what the key ideas are.  We learned this argument from Thomas Jordan and  Natalia Jurga, although it is based on a proof  from \cite{pu} which applies to more general constructions without an underlying grid structure.  \\

\emph{Sketch proof.} Let $r>0$ be very small and $x \in F$ be a `$\mu$-typical point', which is to say we have selected it from a set of large $\mu$ measure guaranteeing it to behave as expected.  This is made precise using the ergodic theorem and Egorov's theorem.  Let $k$ be such that $r \approx n^{-k}$ and let $E \subseteq F$ be the $k$th-level cylinder  (basic rectangle after $k$ steps in the IFS construction) which contains $x$.  Since $x$ is typical, we can assume this is uniquely defined.  We wish to estimate the measure of $B(x,r)$ which, since the height of $E$ is $\approx r$, is roughly the measure of a vertical strip of width $r$ passing through $E$ `centred' at $x$.  Since $x$ is typical, the measure of $E$ is roughly $\exp(-k h(\mu))$ (recall the Shannon-McMillan-Breiman theorem).  Moreover, using the self-affinity of $\mu$, the proportion of the measure of $E$ which lies inside the vertical strip we are interested in will be (roughly) the same as the proportion of an $r m^k$-ball centred at a $\pi \mu$-typical point $x'$ in the projection $\pi F$. The $m^{k}$ factor comes from scaling the base of $E$ up to match the unit interval (where $\pi \mu$ lives).  Since $\pi \mu$ is a self-similar measure satisfying the open set condition, it is well-known and easily shown that it is exact dimensional with dimension given by entropy divided by Lyapunov exponent, that is,
\[
\dim \pi \mu  = \frac{h(\pi\mu)}{\log m}.
\]
Therefore,
\[
\mu(B(x,r)) \approx \mu(E) \cdot \pi\mu (B(x', rm^k)) \approx \exp(-k h(\mu)) \left( rm^k\right)^{ \frac{h(\pi\mu)}{\log m}}.
\]
Then, using $k \approx -\log r/\log n$,
\[
\frac{\log \mu (B(x,r))}{\log r} \approx \frac{h(\pi \mu)}{\log m} + \frac{h(\mu)-h(\pi \mu)}{\log n}
\]
and the error vanishes as $r \to 0$. \hfill \qed
\\

The (Hausdorff) dimension of a measure cannot exceed the Hausdorff dimension of its support and so it is natural to ask if the Hausdorff dimension of a Bedford-McMullen carpet can be realised as the Hausdorff dimension of an invariant measure.  Define $\mu$ by 
\begin{equation} \label{mcmeas}
p_d = N_i^{(\log m/\log n)-1}/m^{\hd F}
\end{equation}
for all $d$ corresponding to rectangles in the $i$th non-empty column.   Note that  $\sum_{d =1}^N p_d = 1$, by Theorem \ref{bedmcformulae}. Then, applying Theorem \ref{lyformula}, we have 
\[
\dim \mu = \hd F.
\]
This self-affine measure is known as the \emph{McMullen measure} and was first used in \cite{mcmullen}.  In fact, the McMullen measure is the unique invariant probability measure of maximal Hausdorff dimension, see \cite{kp-measures}.

\section{Dimensions of  projections onto lines}

How fractal sets and measures behave under projection onto subspaces is a well-studied and central question in fractal geometry and geometric measure theory.  We refer the reader to \cite{FalconerFraserJin, MattilaSurvey} for an overview of the dimension theory of projections in general.  Here we focus only on the planar case and projections of Bedford-McMullen carpets.  Consider the set of lines $L$ passing through the origin in $\mathbb{R}^2$ and write $\pi_L$ for orthogonal projection from $\mathbb{R}^2$ onto the line $L$. In 1954 Marstrand \cite{marstrand} proved that, for Borel $E \subseteq \mathbb{R}^2$,
\begin{equation} \label{mars}
\hd \pi_L E = \min\{\hd E, 1\}
\end{equation}
for almost all $L$.  Here `almost all' is with respect to the natural length measure on the space of lines.  This result stimulated much further work in the area and recently there has been a lot of activity concerning the exceptional set of lines in \eqref{mars};  especially when the exceptional set can be shown to be empty or small in some other sense, see \cite{pablosurvey}.  Ferguson, Jordan and Shmerkin \cite{fjs} proved the following projection theorem for Bedford-McMullen carpets.
\begin{thm} \label{bmproj}
Let $F$ be a Bedford-McMullen carpet with $\log m/ \log n \notin \mathbb{Q}$.  Then
\[
\hd \pi_L F = \min\{\hd F, 1\}
\]
for \textbf{all} $L$ apart from possibly when $L$ is one of the two principle coordinate axes.  
\end{thm}

The strategy for proving this result is to first establish it in the uniform fibres case and then upgrade to the general case by approximating the carpet from within by subsystems with uniform fibres.   Proving the result in the uniform fibres setting is far from straightforward, but we will focus on the `approximating from within' part of the proof.  This trick seems to have useful applications elsewhere.  
\begin{lma}\label{fjslemma}
Let $F$ be a Bedford-McMullen carpet  and $\eps>0$.  Then  there exists a subset $E \subseteq F$ which   is itself  a Bedford-McMullen carpet with  uniform fibres and, moreover,  satisfies 
\[
\hd E  \geq \hd F -\eps.
\]
\end{lma}
\emph{Sketch proof.}  
Let
\[
p_d = N_i^{(\log m/\log n)-1}/m^{\hd F}
\]
be the probability weights defining the McMullen measure, see \eqref{mcmeas}, that is the unique invariant probability measure of maximal dimension.   For a large integer $k \geq 1$, let 
\begin{equation} 
l(k)=\sum_{d =1}^N \lfloor  k p_d \rfloor 
\end{equation}
and
\begin{eqnarray*}
&\,& \hspace{-3mm}\mathcal{I}(k) =  \Big\{  (d_1, \dots, d_{l(k)}) \in \{1, \dots, N\}^{l(k)} \\
&\,& \hspace{3cm} \ : \ \text{for all $1 \leq d \leq N$, }\#\{t :d_t=d\}= \lfloor  k p_d\rfloor \Big\}.
\end{eqnarray*}
Here, $\mathcal{I}(k)$ is a subset of the $l(k)$th iteration of the defining IFS for $F$ chosen such that   digits appear with the `correct frequency', as determined by the McMullen measure.

Consider the IFS consisting of compositions of maps according to $\mathcal{I}(k)$  and denote its attractor by $E(k)$.  It follows that  $E(k)$ is a Bedford-McMullen carpet with  uniform fibres.  Then Theorem \ref{bedmcformulae} gives
\[
  \hd E(k) = \frac{\log M(k)}{\log m(k)} \, + \, \frac{\log ( N(k)/M(k))}{\log n(k)}
\]
where $m(k) = m^{l(k)}$ and $n(k)=n^{l(k)}$ are the integers defining the grid associated with $E(k)$ and $M(k)$ and $N(k)$ are the number of non-empty columns and the total number of rectangles in the construction of $E(k)$, respectively.  Moreover,
\[
N(k) = \frac{l(k)!}{\prod_{d=1}^N \lfloor  k p_d\rfloor !}  
\]
and
\[
M(k) = \frac{l(k)!}{\prod_{i=1}^{M}\left(\sum_{d: \pi d = i} \lfloor  k p_{d}\rfloor \right)!}.
\]
We can apply Stirling's approximation to estimate the dimension of $E(k)$ from below.  A lengthy but straightforward calculation yields
\begin{eqnarray*}
\hd E(k) &=& \frac{\log l(k)! - \sum_{i=1}^{M}\log \left(\sum_{d: \pi d = i} \lfloor  k p_{d}\rfloor \right)!}{l(k)\log m} \\ \\
&\,& \hspace{15mm} \, + \, \frac{\sum_{i=1}^{M}\log \left(\sum_{d: \pi d = i} \lfloor  k p_{d}\rfloor \right)! -  \sum_{d =1}^N\log \lfloor  k p_d\rfloor !}{l(k)\log n} \\ \\
&\geq & s - \eps(k)
\end{eqnarray*}
where $\eps(k) \to 0$ as $k \to \infty$, proving the lemma.   \hfill \qed \\

We give one further example where Lemma \ref{fjslemma} can be applied.    The \emph{modified lower dimension}, denoted by $\mld$, is a modification of the lower dimension to make it monotone, see \cite{fraserbook}.  Specifically, it is defined as $\mld F = \sup\{ \ld E : E \subseteq F\}$ and for compact sets $F$  we have $\ld F \leq \mld F \leq \hd F$.  Perhaps surprisingly it turns out to be equal to the \emph{Hausdorff} dimension, not the lower dimension, in the case of Bedford-McMullen carpets.  This was first observed in \cite{fyu2}. 
\begin{cor} \label{bedmcpsec2}
Let $F$ be a Bedford-McMullen carpet.  Then
\[
\mld F = \hd F =  \frac{\log \sum_{i=1}^{M} N_i^{\log m/ \log n}}{\log m} .
\]
\end{cor}
Corollary \ref{bedmcpsec2} follows immediately from Lemma \ref{fjslemma}.    This   shows how to construct subsets of carpets which have lower dimension arbitrarily close to the Hausdorff dimension of the original carpet, using the fact that the carpets $E$ provided by Lemma \ref{fjslemma} have uniform fibres. This proves the lower bound and the upper bound always holds.  

Ferguson, Fraser and Sahlsten \cite{ffs} considered the projections of measures supported on Bedford-McMullen carpets.  Here the subsystem argument is not applicable and the proof relies on CP-chains and  theory developed by Hochman and Shmerkin \cite{hochmanshmerkin}.
\begin{thm} \label{bmproj2}
Let $\mu$ be a self-affine measure supported on a Bedford-McMullen carpet with $\log m/ \log n \notin \mathbb{Q}$.  Then
\[
\hd \pi_L \mu = \min\{\hd \mu, 1\}
\]
for \textbf{all} $L$ apart from possibly when $L$ is one of the two principle coordinate axes.  
\end{thm}

Theorem \ref{bmproj2} was subsequently generalised by Almarza \cite{almarza} to include Gibbs measures on  transitive subshifts of finite type.

Related to the dimension theory of orthogonal projections, is the dimension theory of slices, see \cite[Chapter 10]{mattila}.  In the plane, a  \emph{slice}  is the intersection of a given set with a line.  Write  $L^\perp$ for the orthogonal complement of a line $L \subseteq \mathbb{R}^2$ and let $E  \subseteq \mathbb{R}^2$ be a  Borel set with $\mathcal{H}^s(E)>0$.   Then  `in typical directions, there are many big slices', that is, 
\[
\hd E \cap (L^{\perp}+x) \geq s-1
\]
for almost all $L$ and positively many $x \in L$.  Here `positively many' $x$ means there is a set of $x$ of positive length.   Moreover, for an arbitrary Borel set $E \subseteq \mathbb{R}^2$ `in all directions typical  slices cannot be too big', that is,
\[
\hd E \cap (L^{\perp}+x)  \leq \max\{\hd E - 1, 0\}
\]
for all $L$ and almost all $x \in L$. Again, there has been a lot of interest in the dimension theory of slices in specific situations, such as for dynamically defined sets $E$.  Slices of Bedford-McMullen carpets were considered by Algom \cite{algom} and the following was proved.  
\begin{thm} \label{algomthm}
Let $F$ be a Bedford-McMullen carpet with $\log m/ \log n \notin \mathbb{Q}$. Then
\[
\bd F \cap (L^{\perp}+x) \leq \max \{ \ad F  - 1, 0\}
\]
for all $L$ and \textbf{all} $x \in L$ provided $L$ is not one of the two principle coordinate axes.
\end{thm}

\section{Survivor sets, hitting problems, and Diophantine approximation}

In this section we briefly touch upon a large and varied literature concerning the study of dynamically or number theoretically defined subsets of a given fractal.  The literature goes far beyond the Bedford-McMullen setting but, as usual, Bedford-McMullen carpets provide an excellent testing ground for the theory.  For example, the `survivor set problem' studies points which do not fall into a given `hole' under iteration of the $(\times m, \times n)$ dynamics.  The `hitting target problem' considers the complementary phenomenon where one focuses on points which hit a given target.  Then there are various related problems in Diophantine approximation, where one is interested in how well points may be approximated by rationals.  This may be interpreted as hitting a prescribed target.

Ferguson, Jordan and Rams \cite{fjr} considered the survivor set problem as follows.  Let $T$ denote the $(\times m, \times n)$ endomorphism which leaves a given Bedford-McMullen carpet $F$ invariant.  Let $U \subseteq F$ be an open set (in the subspace topology)  and define the \emph{survivor set} by 
\[
F_U = \{ x \in F : T^k(x) \notin U \text{ for all } k\}.
\]
The  dimensions of $F_U$ were considered in  \cite{fjr} where $U$ is a fixed finite collection of open cylinders, or a  shrinking metric ball.  They found that the box dimension is related to the escape rate of the measure of maximal entropy through the hole, and the Hausdorff dimension is related to the escape rate of the measure of maximal dimension.

B\'ar\'any and Rams \cite{baranyrams} considered the shrinking target problem as follows.  Let $B_k \subseteq F$ be a sequence of targets where each $B_k$ is a `dynamically defined rectangle'.  Let
\[
\Gamma = \{ x \in F : T^k(x) \in B_k \text{ for infinitely many } k\}
\]
that is, the set of points which, upon iteration of $T$, hit the (moving) target infinitely often.   The Hausdorff dimension of $\Gamma$ is given in \cite{baranyrams} in terms of various complicated entropy functions.

A point $x \in \mathbb{R}^d$ is said to be \emph{badly approximable} if there exists $c>0$ such that, for all $\textbf{p} = (p_1, \dots, p_d) \in \mathbb{Z}^d$ and all $q \in \mathbb{N}$,
\[
\|x-\textbf{p}/q\|_\infty \geq \frac{c}{q^{1+1/d}}.
\]
That is, the badly approximable numbers  are those for which Dirichlet's theorem can be improved by at most a constant factor.  The badly approximable numbers $B(d)$ have zero Lebesgue measure but full Hausdorff dimension in $\mathbb{R}^d$ and a well-studied problem is to determine how $B(d)$ intersects a given fractal set.  Das, Fishman, Simmons and Urba\'nski \cite{das} proved that 
\begin{equation} \label{badc}
\hd B(2) \cap F = \hd F
\end{equation}
when $F$ is a Bedford-McMullen carpet with at least two non-empty columns and at least two non-empty rows.  Interestingly, there is a connection with the (modified) lower dimension here.  It is proved in  \cite{das}  that for an arbitrary  closed set  $E \subseteq \mathbb{R}^d$ satisfying a natural non-degeneracy condition called `hyperplane diffuseness',  we have 
\[
\hd B(d) \cap E \geq  \ld E.
\]
The result \eqref{badc} may then be deduced from this and Corollary \ref{bedmcpsec2}.

\section{Multifractal analysis}

We saw in Theorem \ref{lyformula} that self-affine measures on Bedford-McMullen carpets are exact dimensional, meaning that the local dimension exists and takes a common value at almost every point.  It is an interesting and difficult problem to study the exceptional set, that is, the set of points where the local dimension is not as expected.  This is a $\mu$-null set, but turns out to have full Hausdorff dimension, $\hd F$.  This type of problem is common in   multifractal analysis, see \cite[Chapter 17]{falconer}.  Let $\alpha \geq 0$ and form multifractal decomposition sets
\[
\Delta(\alpha) = \{ x \in F : \dim_{\textup{loc}} (\mu,x) = \alpha\},
\]
where $\dim_{\textup{loc}} (\mu,x)$ is the local dimension of $\mu$ at $x$, recall \eqref{ldd}. In order to understand the fractal complexity of $\Delta(\alpha)$, and thus $\mu$, define the \emph{Hausdorff} and \emph{packing multifractal spectra} as
\[
f_{\textup{H}}^\mu(\alpha) = \hd \Delta(\alpha)
\]
and
\[
f_{\textup{P}}^\mu(\alpha) = \pd \Delta(\alpha),
\]
respectively. It is not useful to define box or Assouad multifractal spectra here since the sets $\Delta(\alpha)$ tend to be dense in $F$ and therefore we need a dimension which is not stable under taking closure to distinguish between different $\alpha$.

Multifractal analysis has been considered in great detail in the context of self-affine measures on Bedford-McMullen carpets.  These measures constitute one of the most complicated examples where the Hausdorff multifractal spectrum is known and given by an explicit formula. That said, many interesting questions remain.

Closely connected to multifractal analysis is the study of the $L^q$-spectrum.  Given a Borel probability measure $\mu$, the  $L^q$-spectrum of $\mu$ is a function $\tau_\mu: \mathbb{R} \to \mathbb{R}$ which captures the coarse structure of the measure by considering $q$th-moment type expressions.   Many interesting fractal features may be analysed via this function. For example, $\tau_\mu(0)$ coincides with the box dimension of the support of the measure and, provided $\tau_\mu$ is differentiable at $q=1$,  $\hd \mu = -\tau_\mu'(1)$, see \cite{ngai}.  More importantly for us, one always has
\begin{equation} \label{mform}
f_{\textup{H}}^\mu(\alpha) \leq f_{\textup{P}}^\mu(\alpha) \leq \tau_\mu^*(\alpha)
\end{equation}
where $\tau_\mu^*(\alpha)$ is the Legendre transform of $\tau_\mu$.  In many cases of interest there is equality throughout in \eqref{mform} in which case we say the multifractal formalism holds.  For example, this holds for self-similar measures satisfying the open set condition.  Self-affine measures on Bedford-McMullen carpets fail to satisfy the multifractal formalism in this sense  in general but, nevertheless, the Hausdorff multifractal spectrum is known and is given by the Legendre transform of an auxiliary moment scaling function $\beta: \mathbb{R} \to \mathbb{R}$.  This was proved by King \cite{king} assuming an additional separation condition known as the \emph{very strong separation condition}, and in full generality by Jordan and Rams \cite{jordanrams}.  Moreover, the $L^q$-spectrum is also known and given by an explicit formula.  This is due to Olsen \cite{olsen}.
\begin{thm}
Let $\mu$ be a self-affine measure on a Bedford-McMullen carpet. There is an explicitly defined real analytic moment scaling function $\beta: \mathbb{R} \to \mathbb{R}$ such that $f_{\textup{H}}^\mu(\alpha) =\beta^*(\alpha)$.    Moreover, the $L^q$-spectrum is real analytic and given by an explicit formula.  In general $\beta$ and $\tau_\mu$ do not coincide.
\end{thm}
The above theorem provides explicit upper and lower bounds for the packing multifractal spectrum.  The problem of computing the packing multifractal spectrum in general was considered in detail by Reeve \cite{reeve} and Jordan and Rams \cite{jordanrams2} and it turns out to be a subtle problem.    Reeve \cite{reeve} considered multifractal analysis of Birkhoff averages, which is different to the multifractal analysis of local dimensions we consider here.  However, in certain cases they can be related and it was shown in \cite{reeve} that the upper bound given by the Legendre transform of the $L^q$-spectrum is generally not sharp.  In \cite{jordanrams2} it was shown that usually the packing multifractal spectrum does not peak at the packing dimension of the Bedford-McMullen carpet.  This is in stark contrast to the Hausdorff case where the Hausdorff multifractal spectrum always peaks at the Hausdorff dimension of the carpet.  In \cite{jordanrams2} the packing multifractal spectrum was computed for a special  family of self-affine measures supported on Bedford-McMullen carpets.  The carpet was allowed only two non-empty columns and the same Bernoulli weight was associated to each rectangle in the same column.  Within this class they were able to show that the packing multifractal spectrum can be discontinuous as a function of the Bernoulli weights.  Again, this is in stark contrast to the Hausdorff case.

Related to multifractal analysis is the study of the quantisation dimensions of a measure.   These were computed in \cite{marc} for self-affine measures on Bedford-McMullen carpets.  Roughly speaking the problem is to determine how well a measure can be approximated by a collection of point masses (quantised).

\section{Open problems}

We conclude this survey article by collecting some open problems relating to Bedford-McMullen carpets. The following question was explicitly asked  in \cite{falconersurvey, FalconerFraserKempton, frasersurvey, istvan} and seems to be technically challenging.

\begin{ques}
Find a precise formula for the intermediate dimensions $\dim_\theta F$ for $F$ a Bedford-McMullen carpet with non-uniform fibres.
\end{ques}

Theorem \ref{bmproj} completely describes the Hausdorff dimensions of the projections of Bedford-McMullen carpets onto lines provided $\log m /\log n \notin \mathbb{Q}$.  The `rational case' remains open, where it seems unlikely that the conclusion of Theorem \ref{bmproj} holds in general.

\begin{ques}
What can be said about the Hausdorff dimensions of the projections of a Bedford-McMullen carpet onto lines when $\log m/ \log n \in \mathbb{Q}$? What about projections of associated self-affine measures?
\end{ques}

There are many interesting open problems in dimension theory and ergodic theory due to Furstenberg which ask about the independence of $\times 2$ and $\times 3$ actions.  Many of these are formulated in  terms of projections or slices of products of $\times 2$ and $\times 3$ invariant sets, see recent breakthroughs \cite{wuf,pablof}.  Bedford-McMullen carpets (and more general $(\times m, \times n)$ invariant sets) therefore provide a natural extension of many of these conjectures since being $(\times m, \times n)$ invariant is more general than being the product of a $\times m$ invariant set and a $\times n$ invariant set.  Theorems \ref{bmproj}, \ref{bmproj2} and \ref{algomthm} are all examples of this in action.  Many questions and conjectures can be formulated and we highlight one example, implicit in \cite{algom}.

\begin{ques}
Let $F$ be a Bedford-McMullen carpet with $\log m/ \log n \notin \mathbb{Q}$.  Is it true that
\[
\hd (F \cap L) \leq \max \{ \hd F -1, 0\}
\]
for all lines $L$ which are not parallel to the coordinate axes?
\end{ques}

It remains an interesting and challenging open problem to fully describe the multifractal analysis of self-affine measures on Bedford-McMullen carpets in the setting of packing dimension.  The following question was explicitly asked in \cite{olsen, reeve}  and shown to be rather subtle in \cite{jordanrams2}.

\begin{ques}
Find a precise formula for the packing multifractal spectrum   $f_\textup{P}^\mu(\alpha)$ for $\mu$ a self-affine measure on a Bedford-McMullen carpet with non-uniform fibres.
\end{ques}

A compact set $E \subseteq \mathbb{R}^2$ is called \emph{tube null} if it can be covered by a collection of tubes of arbitrarily small total area.  A \emph{tube} is an $\varepsilon$-neighbourhood of a line segment.   If   $\hd E < 1$, then $E$ is immediately tube null since one can find a line $L$ such that the $\mathcal{H}^1(\pi_L(E)) = 0$ and then the collection of tubes can be taken transversal to $L$.  In general, the tubes need not be all in the same direction.  In \cite{tubes} it was shown that Sierpi\'nski carpets $E \subseteq \mathbb{R}^2$   are tube null   provided $\hd E < 2$.  Sierpi\'nski carpets are  constructed in the same way as  Bedford-McMullen carpets but with $m=n$ and, as such, are self-similar rather than (strictly) self-affine.   

The approach in \cite{tubes} does not work for Bedford-McMullen carpets in general.  The non-trivial case is when there are no empty columns and no empty rows.  It seems especially difficult to prove tube-nullity in this case if $\log m /\log n \notin \mathbb{Q}$, since then there are no `special projections', see Theorem \ref{bmproj}.

\begin{ques}
Are Bedford-McMullen carpets  tube null, provided they are not the whole unit square?
\end{ques}

\section*{Acknowledgements}

We thank  Natalia Jurga and Istvan  Kolossv\'ary for making several helpful comments and suggestions. We are also grateful to Amir Algom and Meng Wu for interesting discussions.


\begin{thebibliography}{99}

\bibitem[A20]{algom}
A. Algom.
Slicing theorems and rigidity phenomena for self affine carpets, \emph{Proc. London Math. Soc.},  {\bf 121},  (2020), 312--353.


\bibitem[AH19]{algomhochman}
A. Algom and M. Hochman.
Self embeddings of Bedford-McMullen carpets, 
\emph{Ergodic Theory Dynam. Systems}, {\bf 39}, (2019), 577--603.



\bibitem[A17]{almarza}
J. I. Almarza.
CP-chains and dimension preservation for projections of $(\times m, \times n)$-invariant Gibbs measures,
\emph{Adv. Math.},  {\bf 304}, (2017), 227--265.

\bibitem[BK13]{bandtkaenmaki}
 C. Bandt and A. K\"aenm\"aki.
 Local structure of self-affine sets,
\emph{Ergodic Th.  Dynam. Syst.}, {\bf 33}, (2013), 1326--1337.


\bibitem[B07]{baranski}
K.~Bara\'nski.
Hausdorff dimension of the limit sets of some planar geometric constructions,
{\em Adv. Math.}, {\bf 210}, (2007), 215--245.

\bibitem[B15]{barany}
B. B\'ar\'any.
 On the Ledrappier-Young formula for self-affine measures,
\emph{Math. Proc. Cambridge Philos. Soc.}, {\bf 159}, (2015), 405--432.


\bibitem[BHR19]{baranyinvent}
B. B\'ar\'any, M. Hochman and A. Rapaport. 
Hausdorff dimension of planar self-affine sets and measures,
\emph{ Invent. Math.}, {\bf 216}, (2019), 601--659.


\bibitem[BK17]{baranykaenmaki}
B. B\'ar\'any and A. K\"aenm\"aki.
 Ledrappier-Young formula and exact dimensionality of self-affine measures,
\emph{Adv. Math.}, {\bf 318}, (2017), 88--129.

\bibitem[BKR19]{rossi}
B. B\'ar\'any, A. K\"aenm\"aki and E. Rossi.
Assouad dimension of planar self-affine sets,
\emph{Trans. Amer. Math. Soc.}, (to appear), available at: https://arxiv.org/abs/1906.11007

\bibitem[BK17]{baranyrams}
B. B\'ar\'any and M. Rams.
Shrinking targets on Bedford-McMullen carpets,
{\it Proc. Lond. Math. Soc.}, {\bf 117}, (2018),  951--995.

\bibitem[B84]{bedford}
T. Bedford.
 \emph{Crinkly curves, Markov partitions and box dimensions in self-similar sets},
 PhD thesis, University of Warwick, (1984).

\bibitem[B20]{burrell1}
S. A. Burrell.
Dimensions of fractional Brownian images, preprint, available at: https://arxiv.org/abs/2002.03659


\bibitem[BFF19]{burrell2}
S. A. Burrell, K. J. Falconer and J. M. Fraser.
Projection theorems for intermediate dimensions, 
\emph{J. Fractal Geom.}, (to appear), available at: https://arxiv.org/abs/1907.07632

\bibitem[DFSU19]{das}
T. Das, L. Fishman, D. Simmons and M. Urba\'nski.
Badly approximable points on self-affine sponges and the lower Assouad dimension,
\emph{Ergodic Th.  Dynam. Syst.}, {\bf 39}, (2019),  638--657.





\bibitem[F13]{affinesurvey}
K.~J. Falconer.
Dimensions of Self-affine Sets - A Survey, \emph{Further Developments in Fractals and Related Fields}, Birkh\"auser, Boston, (2013), 115--134.

\bibitem[F14a]{falconer}
K. J. Falconer.
{\em Fractal Geometry: Mathematical Foundations and Applications},
 John Wiley \& Sons, Hoboken, NJ, 3rd. ed., 2014.


\bibitem[F20a]{falconersurvey}
K.~J. Falconer.
Intermediate dimensions - a survey,
preprint, (2020).

\bibitem[FFJ15]{FalconerFraserJin}
K.~J. Falconer, J. M. Fraser and X. Jin.
Sixty Years of Fractal Projections,
\emph{Fractal geometry and stochastics V, (Eds. C. Bandt, K.~J. Falconer  and M. Z\"ahle)}, Birkh\"auser, Progress in Probability, (2015).

\bibitem[FFK19]{FalconerFraserKempton}
K. J. Falconer, J. M. Fraser and T. Kempton.
 Intermediate dimensions, \emph{Math. Z.}, 1432--1823, (2019).

\bibitem[F19a]{feng}
D.-J. Feng.
Dimension of invariant measures for affine iterated function systems, preprint.

\bibitem[FW12]{fengwang}
D.-J. Feng and Y. Wang. A class of self-affine sets and self-affine measures,
\emph{J. Fourier Anal. Appl.}, {\bf 11}, (2005), 107--124.

\bibitem[FFS15]{ffs}
A. Ferguson, J. M. Fraser and T. Sahlsten. 
Scaling scenery of $(\times m, \times n)$ invariant measures,
\emph{Adv. Math.}, {\bf 268}, (2015), 564--602.

\bibitem[FJR15]{fjr}
A. Ferguson,  T. Jordan and M. Rams.
Dimension of self-affine sets with holes,
\emph{Ann. Acad. Sci. Fenn. Math.}, {\bf 40}, (2015), 63--88.


\bibitem[FJS10]{fjs}
A. Ferguson, T. Jordan and P. Shmerkin.
The Hausdorff dimension of the projections of self-affine carpets,
\emph{Fund. Math.}, {\bf 209}, (2010), 193--213.

\bibitem[F12]{boxlike}
J. M. Fraser. On the packing dimension of box-like self-affine sets in
the plane, \emph{Nonlinearity}, {\bf 25}, (2012), 2075--2092.

\bibitem[F14b]{fraserassouad}
J.~M. Fraser.
 Assouad type dimensions and homogeneity of fractals,
 {\em Trans. Amer. Math. Soc.}, {\bf 366}, (2014), 6687--6733.

\bibitem[F19b]{frasersurvey}
J. M. Fraser.
Interpolating between dimensions,
\emph{Fractal geometry and stochastics VI}, to appear, Birkh\"auser, Progress in Probability, (2019), available at: https://arxiv.org/abs/1905.11274

\bibitem[F20b]{fraserbook}
J. M. Fraser.
{\em Assouad Dimension and Fractal Geometry},
 Cambridge University Press, Tracts in Mathematics Series, \textbf{222}, (2020).

\bibitem[FH17]{fraserhowroyd}
J. M. Fraser and D. C. Howroyd.
Assouad type dimensions for self-affine sponges,
 \emph{Ann. Acad. Sci. Fenn. Math.}, {\bf 42}, (2017), 149--174.

\bibitem[FY18a]{fyu}
J. M. Fraser and H. Yu.
New dimension spectra: finer information on scaling and homogeneity,
\emph{Adv. Math.}, {\bf 329}, (2018), 273--328.

\bibitem[FY18b]{fyu2}
J. M. Fraser and H. Yu.
Assouad type spectra for some fractal families,
\emph{Indiana Univ.~Math.~J.},  {\bf 67}, (2018), 2005--2043.



\bibitem[GL92]{lalley-gatz}
D. Gatzouras and S.~P. Lalley.
Hausdorff and box dimensions of certain self-affine fractals,
 {\em Indiana Univ. Math. J.}, {\bf 41}, (1992), 533--568.


\bibitem[HS12]{hochmanshmerkin}
M. Hochman and P. Shmerkin.
Local entropy averages and projections of fractal measures,
\emph{Ann. Math.}, {\bf 175}, (2012), 1001--1059.


\bibitem[JR11]{jordanrams}
T. Jordan and M. Rams.
Multifractal analysis for Bedford-McMullen carpets, 
\emph{Math. Proc. of Camb. Phil. Soc.}, {\bf  150},  (2011), 147--156.


\bibitem[JR15]{jordanrams2}
T. Jordan and M. Rams.
Packing spectra for Bernoulli measures supported on Bedford-McMullen carpets, 
\emph{Fund. Math.}, {\bf  229},  (2015), 171--196.




\bibitem[KKR17]{fibred}
A. K{\"a}enm{\"a}ki, H.  Koivusalo and E. Rossi.  
Self-affine sets with fibred tangents,
\emph{Ergodic Th.  Dynam. Syst.}, {\bf 37}, (2017), 1915--1934.

\bibitem[KOR18]{kaenmaki}
  A.~K{\"a}enm{\"a}ki, T.~Ojala, and E.~Rossi.
Rigidity of   quasisymmetric mappings on self-affine carpets, \emph{Int. Math. Res. Not. IMRN}, \textbf{12}, (2018), 3769--3799.





\bibitem[KP96a]{kp}
R. Kenyon and Y. Peres.
Hausdorff dimensions of sofic affine-invariant sets,
\emph{Israel J. Math.}, {\bf 94}, (1996), 157--178.


\bibitem[KP96b]{kp-measures}
R. Kenyon and Y. Peres.
Measures of full dimension on affine-invariant sets,
\emph{Ergodic Th. Dynam. Systems}, {\bf 16}, (1996), 307--323.

\bibitem[KZ16]{marc}
M. Kesseb\"ohmer and S. Zhu.
On the quantization for self-affine measures on Bedford-McMullen carpets, 
{\it Math. Z.}, {\bf  283}, (2016), 39--58.

\bibitem[K95]{king}
J. F. King.
The singularity spectrum for general Sierpi\'nski carpets,
\emph{Adv. Math.}, {\bf 116}, (1995), 1--11.


\bibitem[K20]{istvan}
I. Kolossv\'ary.
On the intermediate dimensions of Bedford-McMullen carpets,
\emph{preprint},  available at: https://arxiv.org/abs/2006.14366

\bibitem[KS19]{istvan2}
I. Kolossv\'ary and K. Simon.
Triangular Gatzouras-Lalley-type planar carpets with overlaps,
\emph{Nonlinearity}, {\bf 32}, (2019), 3294--3341.







\bibitem[LY85a]{led}
F. Ledrappier and L.-S. Young.
 The metric entropy of diffeomorphisms. I. Characterization of measures satisfying
Pesin’s entropy formula,
\emph{ Ann.  Math.}, {\bf  122}, (1985),  509--539.


\bibitem[LY85b]{led2}
F. Ledrappier and L.-S. Young. The metric entropy of diffeomorphisms. II. Relations between entropy, exponents
and dimension, \emph{ Ann.  Math.}, {\bf  122}, (1985), 540--574.

\bibitem[M11]{mackay}
J. M. Mackay.
Assouad dimension of self-affine carpets,
\emph{Conform. Geom. Dyn.} {\bf 15}, (2011), 177--187.

\bibitem[M54]{marstrand}
J.~M. Marstrand.
 Some fundamental geometrical properties of plane sets of fractional dimensions,
{\em Proc. London Math. Soc.(3)}, {\bf 4}, (1954), 257--302.



\bibitem[M95]{mattila}
P. Mattila.
{\em Geometry of sets and measures in Euclidean spaces},
 Cambridge studies in advanced mathematics, \textbf{44}, Cambridge University Press, (1995).


\bibitem[M14]{MattilaSurvey}
P. Mattila.
Recent progress on dimensions of projections,
in {\em Geometry and Analysis of Fractals}, D.-J. Feng and K.-S. Lau (eds.), pp 283--301,
{\it  Springer Proceedings in Mathematics \& Statistics.} {\bf 88}, Springer-Verlag, Berlin Heidelberg, (2014).


\bibitem[M84]{mcmullen}
C. McMullen.
 The Hausdorff dimension of general Sierpi\'nski carpets,
 {\em Nagoya Math. J.}, {\bf 96}, (1984), 1--9.


\bibitem[N97]{ngai}
S.-M. Ngai.
A dimension result arising from the $L^q$-spectrum of a measure,
\emph{Proc. Amer. Math. Soc.}, {\bf 125}, (1997), 2943--2951.

\bibitem[O98]{olsen}
L. Olsen.
Self-affine multifractal Sierpi\'nski sponges in $\mathbb{R}^d$,
{\em Pacific J. Math}, {\bf 183}, (1998), 143--199.



\bibitem[P94a]{perespacking}
Y. Peres.
The packing measure of self-affine carpets,
\emph{Math. Proc. Cambridge Philos. Soc.}, {\bf  115}, (1994), 437--450.

\bibitem[P94b]{pereshausdorff}
Y. Peres.
The self-affine carpets of McMullen and Bedford have infinite Hausdorff measure,
\emph{Math. Proc. Cambridge Philos. Soc. }, {\bf 116},  (1994),  513--526.

\bibitem[PU89]{pu}
F. Przytycki and M.  Urba\'nski.
On the Hausdorff dimension of some fractal sets,
\emph{Studia Math.}, {\bf  93},  (1989),  155--186.

\bibitem[PSSW20]{tubes}
A. Py\"or\"al\"a, P.  Shmerkin, V. Suomala and  Meng Wu.
Covering the Sierpi\'nski carpet with tubes,
preprint, available at: https://arxiv.org/abs/2006.00499.

\bibitem[R12]{reeve}
 H. W. J. Reeve.
The packing spectrum for Birkhoff averages on a self-affine repeller, 
\emph{Ergodic Th. Dynam. Syst.}, {\bf 32}, (2012),  1444--1470.

\bibitem[S15]{pablosurvey}
 P.  Shmerkin.
Projections of self-similar and related fractals: a survey of recent developments.\emph{Fractal geometry and stochastics V, (Eds. C. Bandt, K.~J. Falconer  and M. Z\"ahle)}, Birkh\"auser, Progress in Probability, (2015).

\bibitem[S19]{pablof}
 P.  Shmerkin.
On Furstenberg's intersection conjecture, self-similar measures, and the $L^q$ norms of convolutions, \emph{ Ann. of Math.}, {\bf 189}, (2019),  319--391.

\bibitem[W19]{wuf}
M. Wu.
A proof of Furstenberg's conjecture on the intersections of $\times p$  and $\times q$-invariant sets, \emph{Ann. of Math.}, {\bf 189}, (2019), 707--751.

\end{thebibliography}
\end{document}